\documentclass[12pt]{amsart}

\usepackage{verbatim}

\usepackage{amscd,amssymb,fullpage,epsf}
\usepackage[mathscr]{eucal}
\usepackage{graphicx}

\theoremstyle{plain}

\numberwithin{equation}{section}

\newtheorem{theorem}{Theorem}[section]
\newtheorem{lemma}[theorem]{Lemma}
\newtheorem{proposition}[theorem]{Proposition}
\newtheorem{corollary}[theorem]{Corollary}

\theoremstyle{definition}

\newtheorem{definition}{Definition}[section]

\theoremstyle{remark}

\newtheorem*{example*}{Example}
\newtheorem*{remark*}{Remark}

\newcommand{\R}{\mathbb{R}}
\newcommand{\C}{\mathbb{C}}
\newcommand{\Z}{\mathbb{Z}}
\newcommand{\Chat}{\widehat{\mathbb{C}}}

\newcommand{\p}{\partial}
\newcommand\opart{{\overline{\partial}}}

\newcommand\of{{\overline{f}}}
\newcommand\oh{{\overline{h}}}
\newcommand\oz{{\overline{z}}}
\newcommand\ozeta{{\overline{\zeta}}}
\newcommand\oq{\overline{q}}

\newcommand\otau{\overline{\tau}}
\newcommand\ow{\overline{w}}
\newcommand{\Teich}{\textit{Teich}}
\newcommand{\QF}{\textit{QF}}\newcommand{\M}{\mathcal{M}}
\renewcommand{\H}{\mathbb{H}}
\newcommand{\dbar}{\overline{\partial}}
\DeclareMathOperator{\Spec}{Spec}
\DeclareMathOperator{\Tr}{Tr}
\let\Re=\undefined
\DeclareMathOperator{\Re}{Re}
\let\Im=\undefined
\DeclareMathOperator{\Im}{Im}

\newcommand{\detp}{\operatorname{det}'}
\DeclareMathOperator{\Diff}{\textit{Diff}}
\DeclareMathOperator{\horizontal}{\mathsf{v}}
\DeclareMathOperator{\vertical}{\mathsf{w}}

\newcommand{\ssubset}{\subset\subset}

\newcommand{\omu}{\overline{\mu}}
\newcommand{\onu}{\overline{\nu}}
\newcommand{\fm}{f^{\mu}}

\newcommand{\ofm}{\overline{\fm}}

\newcommand{\fmn}{f_{\mu,\nu}}
\newcommand{\ofnm}{\overline{f_{\nu, \mu}}}
\newcommand{\pfm}{\partial \fm}
\newcommand{\opfm}{\overline{\partial} \fm}

\newcommand{\pfmn}{\p \fmn}

\newcommand{\opfmn}{\opart \fmn}
\newcommand{\opofnm}{\opart \ofnm}

\newcommand{\mat}{ \begin{pmatrix} }
\newcommand{\emat}{\end{pmatrix}}

\newcommand{\op}{{\overline{\partial}}}
\newcommand{\oalpha}{\overline{\alpha}}

\renewcommand{\symbol}{\sigma_2}
\DeclareMathOperator{\dist}{dist}
\renewcommand{\S}{P}
\newcommand{\T}{Q}
\newcommand{\zhat}{\hat{z}}
\newcommand{\omegaWP}{\omega_{\textit{WP}}}

\begin{document}

\title[Holomorphic extensions of Laplacians, Determinants]
{Holomorphic Extensions of Laplacians and Their Determinants}

\author{Young-Heon Kim}

\address{Department of Mathematics\\
  Northwestern University   \\
  Evanston, IL 60208 \\
  USA}

\email{huns@math.northwestern.edu}

\date{\today}

\thanks{This research was partially supported by
  the NSF under grant DMS-9704320}


\begin{abstract}


The Teichm\"uller space $\Teich(S)$ of a surface $S$ in genus $g>1$ is
a totally real submanifold of the quasifuchsian space $\QF(S)$. We
show that the determinant of the Laplacian $\detp(\Delta)$ on
$\Teich(S)$ has a unique holomorphic extension to $\QF(S)$.  To
realize this holomorphic extension as the determinant of differential
operators on $S$, we introduce a holomorphic family $\{
\Delta_{\mu,\nu} \}$ of elliptic second order differential operators
on $S$ whose parameter space is the space of pairs of Beltrami
differentials on $S$ and which naturally extends the Laplace operators
of hyperbolic metrics on $S$.  We study the determinant of this family
$\{ \Delta_{\mu,\nu} \}$ and show how this family realizes the
holomorphic extension of $\detp(\Delta)$ as its determinant.
\end{abstract}

\maketitle


\section{Introduction}

In this paper, we discuss determinants of Laplacians of Riemann
surfaces and their holomorphic extensions.

Given a closed Riemannian manifold $X$ with metric $m$, its
corresponding Laplacian $\Delta$ is a self-adjoint positive definite
elliptic second order differential operator on functions on $X$,
which has discrete spectrum
\begin{equation*}
  \lambda_0 = 0 < \lambda_1 \leq \lambda_2 \leq \cdots \leq
  \lambda_k \leq \cdots \rightarrow \infty .
\end{equation*}
The determinant of the operator $\Delta$ may be defined formally as
the product of the nonzero eigenvalues of $\Delta$. A regularization
$\detp(\Delta)$ of this product was defined by Ray and Singer
\cite{RS} \cite{RS2},  using the zeta function of $\Delta$.

This determinant $\detp(\Delta)$ has appeared to be
very important in mathematics.
For example, in \cite{OPS}, (see also \cite{Sa}), Osgood, Phillips and Sarnak
studied $-\log\detp(\Delta)$ as a ``height'' function on the space of
metrics on a compact orientable smooth surface $S$ of genus $g$. For
$g>1$, they showed that when restricted to a given conformal class of
metrics on $S$, it attains its minimum at the unique hyperbolic metric
in this conformal class, and has no other critical points. Thus, to
find Riemannian metrics on $S$ which are extremal, in the sense that
they minimize $-\log\detp(\Delta)$, it suffices to consider its
restriction to the moduli space $\M_g$ of hyperbolic metrics on a
Riemann surface $S$ of genus $g$.
It was shown by Wolpert that this restriction is a proper function
(see \cite{W4}), which was used also by Osgood, Phillips and Sarnak
to show that the isospectral sets (with respect to the Laplacian)
of isometry classes of metrics  on $S$ are all compact in the
$C^\infty$ topology (see \cite{OPS2}).

The universal cover of the orbifold $\M_g$, with covering group the
mapping class group $\Gamma_g$, is the Teichm\"uller space
$\Teich(S)$. The function $-\log\detp(\Delta)$ lifts to a function on
the Teichm\"uller space $\Teich(S)$ invariant under $\Gamma_g$.  In
the first part of this thesis,
we are interested in the function theoretic properties of
$\log\detp(\Delta)$ on $\Teich(S)$.

\subsection{Holomorphic extensions of determinants of Laplacians}

Before stating our first main theorem, consider the special case of
genus $1$.

\begin{example*}[\cite{RS2} or \cite{S1}, p.\ 33, (A.1.7)]
  For $z \in \H$, let $T_z$ be the flat torus obtained by the
  lattice of $\C$ generated by $1$ and $z$. Then the
  determinant of Laplacian of this flat torus is
  \begin{equation*}
  \log \detp(\Delta)(z) = \log (2\pi (\text{Im} \ z )^{1/2} |\eta (z) |^2 )
  \end{equation*}
  where $\eta(z)= q^{1/24} \prod_{n=1}^{\infty} (1- q^{n})$ for
  $q=e^{2\pi i z}$ is the Dedekind eta function; this is a modular form of
  weight $1/2$.

  The manifold $\H$ has a complexification $\H \times
  \overline{\H}$, and the function $\log\detp(\Delta)(z)$ on
  the diagonal $\{w=\overline{z}\}$ has a unique holomorphic extension to
  $\H \times \overline{\H}$, namely,
  \begin{equation*}
    \log \bigl( 2\pi (\frac{z-w}{2i})^{1/2} \, \eta(z) \,
    \overline{\eta(\overline{w})}\bigr) .
  \end{equation*}
\end{example*}

We show that even in higher genus $g>1$, the function $\log
\detp(\Delta)$ has a unique holomorphic extension. In higher genus,
the objects corresponding to $\H$ and $\H\times\overline{\H}$ are the
Teichm\"uller space $\Teich(S)$ and the quasifuchsian space
\begin{equation*}
  \QF(S)=\Teich(S) \times \Teich(\overline{S}) \cong
  \Teich(S) \times \overline{\Teich(S)},
\end{equation*}
respectively where the real analytic manifold $\Teich(S)$ imbeds as
the diagonal in $\QF(S)$. Bers's ``simultaneous uniformization
theorem'' \cite{B} identifies the quasifuchsian space $\QF(S)$ with
the space of hyperbolic metrics modulo isotopies on the $3$-manifold
$S \times \R$, whose ideal boundary at infinity is conformally
isomorphic to a pair of Riemann surfaces. McMullen recently used the
quasifuchsian space to study the geometry of the Teichm\"uller space
via the above complexification \cite{M}.

Now let us state our first main result.
\begin{theorem} \label{main}
  The function $\log\detp(\Delta)$ on $\Teich(S)$ has a unique
  holomorphic extension to the quasifuchsian space $\QF(S)$.
\end{theorem}

In the proof of Theorem~\ref{main}, we use the Belavin-Knizhnik
formula (see Theorem~\ref{curvature}), proved by Wolpert \cite{W3} and
by Zograf and Takhtajan \cite{ZT} and the holomorphic extension of the
Weil-Petersson form constructed by Platis~\cite{P} (see
Theorem~\ref{Platis}).

We remark that the asymptotic behavior of $\log\detp(\Delta)$ near the
boundary of Teichm\"uller space is important in both geometry and
physics and was studied in \cite{W4} and \cite{BB}. It would be
interesting to understand the asymptotic behavior of the holomorphic
extension of $\log\detp(\Delta)$ near the boundary of the
quasifuchsian space.

In view of Theorem~\ref{main}, it is natural to ask whether there is
an actual family of elliptic differential operators on $S$ whose
determinant realizes the holomorphic extension of $\detp (\Delta)$.
To address this question
we introduce  a family $\{\Delta_{\mu,\nu}
\}$ of elliptic second order differential operators on $S$ which
is holomorphic with respect to its parameter $(\mu,\nu)$,
the pair of Beltrami differentials and which uniquely extends the
Laplacians of hyperbolic metrics.  Because of holomorphy of this
family, the differential operators $\Delta_{\mu,\nu}$ cannot be
self-adjoint off the diagonal $\{ \mu = \nu \}$.  These operators $
\Delta_{\mu,\nu} $ are new examples of non-self-adjoint elliptic
second order differential operators with a natural geometric origin!

\subsection{Holomorphic extensions of Laplacians and their determinants} \label{SS:main theorem}

To state our theorem on the holomorphic extension $\Delta_{\mu , \nu}$
of Laplacians we need a few terminologies.  Recall that a marking on
$S$ is a Riemann surface $X_0$ together with an oriented
diffeomorphism between $X_0$ and $S$. A Beltrami differential $\mu$ on
$X_0$ is a complex $(-1,1)$-form which in one (and hence all) local
representations
\begin{equation*}
  \mu = \mu(z) \frac{d\oz}{dz}
\end{equation*}
satisfies $\|\mu\|_\infty < 1$. The space $M(X_0)$ of smooth Beltrami
differentials on $X_0$ is a contractible complex analytic manifold
modeled on a Fr\'echet space.  Denote by $M(S)$ the space of smooth
complex structures on $S$, which is equivalent by the uniformization
theorem to the space of hyperbolic metric on $S$. Then $M(X_0)$ gives
a complex coordinate chart on $M(S)$, in which the origin $0\in
M(X_0)$ corresponds to $X_0 \in M(S)$ (see \cite{EE}).  Denote the
complex conjugate of $M(X_0)$ by $\overline{M(X_0)}$. The diagonal
\begin{equation*}
  \{ (\mu, \mu) \mid \mu \in M (X_0) \} \subset
  M(X_0)\times\overline{M(X_0)}
\end{equation*}
is a totally real submanifold. Given $0<k<1$ and $E>0$, we introduce the
space of Beltrami differentials
\begin{equation*}
  M_{k,E}(X_0) = \{ \mu \in M(X_0) \mid \text{$\|\mu\|_\infty
    < k$ and $\|\mu\|_{C^2 (X_0 )}<E$}\} .
\end{equation*}
where the $C^2$-norm $\| \cdot \|_{C^2 (X_0 )}$ is defined by the
hyperbolic metric on $X_0$.

The upper-half plane $\H$ with its standard hyperbolic metric $y^{-2}
(dx^2+dy^2)$ is the Riemannian universal cover of $X_0$; the covering
transformation group $G$ is called the Fuchsian group of $X_0$.  The
Laplacian of $\H$ is given by the formula
\begin{equation*}
  \Delta_\H = (z - \oz)^2 \frac{\p^2}{\p z \p \oz} ,
\end{equation*}
where $z$ is the standard coordinate of $\H$, and it induces
the Laplacian $\Delta$ of
the hyperbolic surface $X_0 = \mathbb{H}/G$.

 Denote by $M^G$ the set of Beltrami differentials on $\H$
which transform as
\begin{equation*}
  \mu(z) = \mu(g(z)) \frac{\overline{\p g}}{\p g}
\end{equation*}
for all $g\in G$. Then $M(X_0)$ is identified with $M^G$.  It is well
known that for each Beltrami differential $\mu$ on $\H$ there exists
unique quasiconformal homeomorphism $f^\mu \, : \, \H \to \H$
satisfying the Beltrami differential equation
\begin{equation*}
  \op f = \mu\, \p f
\end{equation*}
whose continuous extension to the real axis fixes $0,1,\infty$.

We are now ready to state our second main theorem.
\begin{theorem} \label{T:main} There exists unique family of elliptic
  second order differential operators $\Delta_{\mu,\nu}$ on $S$
  parametrized by $(\mu,\nu) \in M(X_0)\times \overline{M(X_0)}$, with
  the following properties:
  \begin{enumerate}
  \item $\Delta_{\mu,\nu}$ depends holomorphically on $(\mu,\nu)$;
  \item the lift of $\Delta_{\mu,\mu}$ to $\H$ is the pull-back of the
    Laplacian $\Delta_\H$ by the quasiconformal mapping
    $f^\mu:\H\to\H$, i.e., $\Delta_{\mu , \mu}$ is the
    Laplacian of the hyperbolic
    metric on $S$ induced by the pullback hyperbolic metric on
    $\mathbb{H}$ by the map
    $f^\mu$;
  \item given $0<k<1$ and $E>0$, there exists a constant $\epsilon>0$
    such that if $\mu$, $\nu\in M_{k,E}(X_0)$ and
    \begin{equation*}
      \|\mu-\nu\|_{C^2(X_0)} < \epsilon ,
    \end{equation*}
    the determinant $\detp (\Delta_{\mu, \nu})$  is defined, and depends
    holomorphically on $(\mu,\nu)$.
  \end{enumerate}
\end{theorem}

The operator $\Delta_{\mu,\nu}$ is constructed by modifying the
explicit expression for $(f^\mu)^*\Delta_\H$, incorporating the
quasifuchsian parameter $(\mu,\nu)$ and corresponding quasiconformal
mapping $f_{\mu,\nu}$. We use a result of Ahlfors and Bers \cite{AB},
that the unique normalized solution of Beltrami differential equation
depends analytically on the Beltrami differential.

To establish property (3), we apply the definition of determinant
using complex powers of elliptic operators due to Seeley (\cite{Se1},
\cite{Se2}, \cite{Sh}, and \cite{KV}). The restriction
$\|\mu-\nu\|_{C^2(X_0)}<\epsilon$ is introduced to satisfy the
conditions for the construction of complex power.

Denote by $\widetilde{\detp (\Delta)}$ the holomorphic extension of
$\detp(\Delta)$ to $\QF(S)$ obtained in Theorem~\ref{main}. We have
the principal fiber bundle
\begin{equation} \label{fibre Teich}
  \begin{CD}
    \Diff_0(S) @>>> M(X_0) \\
    @. @VV\pi V  \\
    @. \Teich(S) ,
  \end{CD}
\end{equation}
where the projection $\pi$ is known to be holomorphic (see \cite{EE}).
This gives rise to the principal fiber bundle
\begin{equation} \label{fibre QF}
  \begin{CD}
    \Diff_0 (S) \times \Diff_0 (S) @>>> M(X_0) \times \overline{M(X_0)} \\
    @. @VV{\pi\times\bar{\pi}}V  \\
    @. \QF(S) .
  \end{CD}
\end{equation}
The lift $(\pi\times\bar{\pi})^* \widetilde{\detp (\Delta)}$ is
holomorphic on $M(X_0) \times \overline{ M(X_0)}$. We know by
Theorem~\ref{T:main} (2) that
\begin{equation*}
  \detp (\Delta_{\mu , \mu})= (\pi\times\bar{\pi})^*
  \widetilde{\detp(\Delta)} (\mu , \mu) ,
\end{equation*}
and by Theorem~\ref{T:main} (3) that the determinant $\detp
(\Delta_{\mu , \nu})$ is defined and holomorphic on some open
neighborhood $N$ of the diagonal in $M(X_0) \times \overline{M(X_0)}$.
Therefore, by analytic continuation, we have the equality
\begin{equation*}
  \detp (\Delta_{\mu , \nu}) =
  (\pi\times\bar{\pi})^* \widetilde{\detp (\Delta)} (\mu , \nu) \quad
  \text{for $(\mu , \nu ) \in  N$,}
\end{equation*}
and we may regard the holomorphic function
$(\pi\times\bar{\pi})^*\widetilde{\detp(\Delta)}$ as the determinant
of $\Delta_{\mu,\nu}$ even for those $(\mu,\nu)$ to which
Theorem~\ref{T:main} (3) does not apply. That is, on all of
$M(X_0)\times\overline{M(X_0)}$, we may define
\begin{equation}\label{new def det}
  \detp (\Delta_{\mu , \nu}) =
  (\pi\times\bar{\pi})^*\widetilde{\detp(\Delta)}(\mu,\nu) .
\end{equation}

\begin{remark*}
  From the family $\{ \Delta_{\mu,\nu}\}$, we may construct
  holomorphic families of elliptic operators in a neighborhood of each
  Teichm\"uller point $([X_0],[X_0])$ in $QF(S)$, using the
  Ahlfors-Weill section $s$ of the fibre bundle \eqref{fibre Teich}
  (see \cite{AW} or \cite{IT} pp.\ 153-157).
  This induces a holomorphic section $s \times \bar{s} $ of fibration $\pi\times\bar{\pi}$ of
  \eqref{fibre QF}, defined in a neighborhood $U$ of the point $([X_0], [X_0])$
  in $\QF(S)$.  Clearly, by \eqref{new def det},
  \begin{equation*}
    \detp (\Delta_{(s[X],\bar{s}[Y])}) = \widetilde{\detp (\Delta)}
    ([X], [Y]) \quad \text{on $U$.}
  \end{equation*}
  However, this method does not give rise to a family of operators
  over all of $\QF(S)$, since by Earle \cite{E}, there is no global
  holomorphic cross-section for the fibre bundle
  $\pi:M(X_0)\to\Teich(S)$ of \eqref{fibre Teich}.
\end{remark*}

\subsection*{Plan of the paper}

In Section~\ref{Ch:det}, we prove Theorem~\ref{main} and in
Section~\ref{Ch:Delta}, we prove Theorem~\ref{T:main}. In subsequent
sections, we provide proof of the results used in Section~\ref{Ch:Delta}.

\subsection*{Acknowledgment}
This paper is author's thesis for Northwestern University.
He deeply thanks his thesis advisor Ezra Getzler for guidance and
support. He is also grateful to
Curtis McMullen, Peter Sarnak, Andr\'as Vasy, and Jared Wunsch
for helpful discussions or comments.

\section{Holomorphic Extensions of Determinants of Laplacians}
\label{Ch:det}

In this section, we use several fundamental facts of Teichm\"uller spaces
and the determinant of Laplacians to prove
Theorem~\ref{main}.

\subsection{Preliminaries}
In this subsection, we review the facts that we need on Teichm\"uller
spaces
and quasifuchsian spaces, including the Belavin-Knizhnik formula and
Platis's theorem. In the next subsection, we prove Theorem~\ref{main}.

\subsection*{Determinants of Laplacians}

Let $\Delta $ be the Laplace-Beltrami operator on functions on a
compact Riemannian manifold $M$. Let
\begin{equation} \label{zeta}
  \zeta_\Delta(s)=  \sum_{\lambda \in \Spec(\Delta)\setminus\{0\}}
  \lambda^{-s}
\end{equation}
be the zeta-function of $\Delta$. The determinant $\detp(\Delta)$ is
defined (see \cite{RS}) as
\begin{equation} \label{det'}
  - \log \detp(\Delta) = \frac{d\zeta_\Delta(0)}{ds} .
\end{equation}
The
sum in \eqref{zeta} is absolutely convergent for $\text{Re}\, s > \frac{\dim M}{2}$
sufficiently large, and has a meromorphic extension to the whole
complex plane. This meromorphic extension is regular at $s=0$, and so
there is no difficulty in taking the derivative at $s=0$ in
\eqref{det'}.

\subsection*{Teichm\"uller spaces}
A general reference for this section is \cite{IT}.

Let $S$ be an oriented closed surface with genus $g > 1$. The
Teichm\"uller space $\Teich(S)$ of $S$ is the space of isotopy classes
of hyperbolic Riemannian metrics on $S$, that is, metrics with
Gaussian curvature $-1$. By uniformization theorem,
$\Teich(S)$ is also the space of isotopy classes of complex
structures on $S$.

The set of equivalence classes of hyperbolic metrics (or equivalently
complex structures) under orientation preserving diffeomorphisms on
$S$ forms the moduli space $\M_g$ of compact Riemann surfaces of genus
$g$.

Denote the group of orientation preserving diffeomorphisms on $S$ by
$\Diff^+(S)$, and the group of isotopies by $\Diff_0(S)$. The mapping
class group
\begin{equation*}
  \Gamma_g = \Diff^+(S) / \Diff_0 (S)
\end{equation*}
is a discrete group which acts properly discontinuously on
$\Teich(S)$. Thus $\Teich(S)$ is almost a covering space of $\M_g$, with
covering transformation group $\Gamma_g$:
\begin{equation*}
  \begin{CD}
    \Gamma_g @>>> \Teich(S)  \\
    @. @VVV \\
    @. \M_g = @. \Gamma_g \backslash \Teich(S)
  \end{CD}
\end{equation*}
The only caveat is that the action of $\Gamma_g$ is not free, i.e.\
there are points in $\Teich(S)$ which are fixed under some finite
subgroups of $\Gamma_g$. These points descend to $\M_g$ as orbifold
singularities.

Fixing a hyperbolic metric on $S$, we may decompose $S$ into $2g-2$
pairs of pants, separated by closed geodesics
$\gamma_1,\dots,\gamma_{3g-3}$.
A hyperbolic pair of pants is determined up to isometry by the lengths
of its boundary geodesics. Given the combinatorial pants decomposition
of $S$, we get a hyperbolic metric by specifying the lengths $l_i$
($l_i >0 $) of the geodesics $\gamma_i$ and the angle $\theta_i$ by
which they are twisted along $\gamma_i$ before gluing. Let
$\tau_i=l_i\theta_i/2\pi$, $i=1,\dots,3g-3$. Then the system of
variables
\begin{equation*}
  (l_1 , \dots , l_{3g-3}, \tau_1 , \dots , \tau_{3g-3} )
\end{equation*}
is a real analytic coordinate system on $\Teich(S)$, called the
\emph{Fenchel-Nielsen coordinates} of $\Teich(S)$. This coordinate
system gives a diffeomorphism
\begin{equation*}
  \Teich(S) \approx \mathbb{R}_+ ^{3g-3} \times \mathbb{R}^{3g-3} \, .
\end{equation*}

There is a a natural symplectic form $\omegaWP$ on $\Teich(S)$,
called the Weil-Petersson form. By a theorem of Wolpert (\cite{W1},
\cite{W2}; see also \cite{IT}), this form is given in Fenchel-Nielsen
coordinates by the formula
\begin{equation} \label{WP}
  \omegaWP = \sum_{i=1}^{3g-3} dl_i \wedge \, d\tau_i .
\end{equation}

The Teichm\"uller space $\Teich(S)$ has a natural complex structure,
for which $\omegaWP$ is a K\"ahler form. The following theorem is
well known. (See, for example, \cite{A}.)
\begin{theorem} \label{Teich as domain}
  For a closed surface $S$ with genus $g>1$, $\Teich(S)$ is
  biholomorphic to a bounded open contractible domain in $\C^{3g -3}$.
\end{theorem}

\begin{corollary} \label{global coordinate}
  There are global holomorphic coordinates $z=(z^1,\dots,z^{3g-3})$ on
  $\Teich(S)$.
\end{corollary}

\subsection*{Quasifuchsian spaces}
While Teichm\"uller space is a space of Riemann surfaces, the
quasifuchsian space defined by Lipman Bers (See \cite{B}) is a space
of pairs of Riemann surfaces. The quasifuchsian space $\QF(S)$ of the
surface $S$ may simply be defined as
\begin{equation*}
  \QF(S) = \Teich(S) \times \Teich(\overline{S}) .
\end{equation*}
Here, $\overline{S}$ denotes the surface $S$ with the opposite
orientation.

The complex conjugate $\overline{X}$ of a Riemann surface $X$ is
defined by the following diagram:
\begin{equation} \label{conjugation}
  \begin{array}{clcl}
    \H & \longrightarrow & \overline{\H} & \\
    \downarrow & & \downarrow & \\
    X & \dashrightarrow & \overline{X} & \\
  \end{array}
\end{equation}
The upper arrow is complex conjugation, and the vertical arrows are
the universal coverings given by the uniformization theorem for
Riemann surfaces. There is a canonical map from $\Teich(S)$ to
$\Teich(\overline{S})$ defined by sending a Riemann surface
$X\in\Teich(S)$ to its complex conjugate
$\overline{X}\in\Teich(\overline{S})$.  As complex manifolds,
$\Teich(\overline{S}) \cong \overline{\Teich(S)}$, where
$\overline{\Teich(S)}$ is the complex conjugate of $\Teich(S)$, i.e.\
the holomorphic structure of $\overline{\Teich(S)}$ is the
anti-holomorphic structure of $\Teich(S)$.

The diagonal map
\begin{equation*}
  \Teich(S) \hookrightarrow \Teich(S) \times \overline{\Teich(S)}
\end{equation*}
sending $X\in\Teich(S)$ to $(X,X)$ embeds $\Teich(S)$ as a totally
real submanifold into $\QF(S)$. The action of $\Gamma_g$ on
$\Teich(S)$ extends to $\QF(S)$ by the diagonal action: for $\rho \in
\Gamma_g$ and $(X,Y) \in \QF(S)=\Teich(S)\times\overline{\Teich(S)}$,
\begin{equation*}
  \rho\cdot(X, Y) = (\rho\cdot X , \rho \cdot Y) .
\end{equation*}

By Corollary~\ref{global coordinate}, $\QF(S)=\Teich(S) \times
\overline{\Teich(S)}$ has global holomorphic coordinates
\begin{equation*}
  (z^1 , \dots , z^{3g-3} , w^1 , \dots , w^{3g-3}) .
\end{equation*}
We abbreviate this coordinate system to $(z,w)$.
Then $\Teich(S)=\{w=\overline{z}\}\subset\QF(S)$.

\subsection*{Holomorphic extension of Weil-Petersson form}

The following result is due to Platis (\cite{P}, Theorems 6 and 8).
\begin{theorem} \label{Platis}
  The differential form $i\omegaWP$ on the Teichm\"uller space
  $\Teich(S)$ has an extension $\Omega$ to the quasifuchsian space
  $\QF(S)$ which is a holomorphic non-degenerate closed $(2,0)$-form
  whose restriction to the diagonal
  $\Teich(S)\subset\QF(S)\cong\Teich(S)\times\overline{\Teich(S)}$ is
  $i\omegaWP$.
\end{theorem}

The following lemma is elementary.
\begin{lemma} \label{identity theorem}
  Let $U\subset\C^n$ be a connected complex domain, and let $\phi$ be
  a holomorphic function on $U\times\overline{U}$ whose restriction to
  the diagonal $U\subset U\times\overline{U}$ vanishes. Then $\phi$
  vanishes on all of $U\times\overline{U}$.
\end{lemma}

We can now prove the following result.
\begin{proposition}
  \label{form of Omega}
  In terms of the holomorphic coordinate system
  \begin{equation*}
    (z,w)=(z^1, \dots, z^{3g-3}, w^1, \dots , w^{3g-3} )
  \end{equation*}
  on $\Teich(S) \times \overline{\Teich(S)}$, the $2$-form $\Omega$ of
  Theorem~\ref{Platis} may be written locally as
  \begin{equation*}
    \Omega = \sum_{i,j} \Omega_{ij} \, dz^i \wedge dw^j .
  \end{equation*}
\end{proposition}
\begin{proof}
  Since $\Omega$ is $(2,0)$ form, we may write
  \begin{equation*}
    \Omega = \sum_{i,j} \bigl( A_{ij} \, dz^i \wedge dz^j + B_{ij} \,
    dz^i \wedge dw^j + C_{ij} \, dw^i \wedge dw^j \bigr) .
  \end{equation*}
  Because the restriction $i\omegaWP$ of $\Omega$ to the diagonal
  $\{w=\oz\}$ is $(1,1)$-form, we see that $A_{ij}$ and $C_{ij}$
  vanish on the diagonal. Since $\Omega$ is holomorphic,
  Lemma~\ref{identity theorem} shows that $A_{ij}$ and $C_{ij}$
  vanish.
\end{proof}

\subsection*{The Laplacian on hyperbolic surfaces and the
  Belavin-Knizhnik formula}

Let $X$ be a compact hyperbolic surface of genus $g>1$, and let
$\Delta$ be the Laplacian on scalar functions on $X$. On the universal
cover $\H$ of $X$, the pull-back of $\Delta$ by the covering map may
be written as
\begin{equation*}
  \Delta = (z-\oz)^2 \frac{\p^2}{\p z\p\oz} .
\end{equation*}

The \emph{Siegel upper half space} $\S_g$ is the space of complex
symmetric matrices in $\C^{g\times g}$ with positive definite
imaginary part. The period matrix $\tau$ is a holomorphic map from
$\Teich(S)$ to $\S_g$.

We will use the Belavin-Knizhnik formula, proved by Wolpert and by
Zograf and Takhtajan. (See \cite{W3} and \cite{ZT}.)  We only need the
following special case of this theorem (\cite{ZT}, Theorem 2).
\begin{theorem} \label{curvature}
  In $\Teich(S)$,
  \begin{equation*}
    \p \dbar \, \log \left( \frac{\detp(\Delta)}{\det(\Im\tau)}
    \right) = - \frac{i}{6 \pi}  \, \omegaWP ,
  \end{equation*}
  where $\Im\tau$ is the imaginary part of the period matrix $\tau$. The
  differential operator $\p\dbar$ comes from the complex structure on
  $\Teich(S)$.
\end{theorem}

This formula and the result of the next section together with the
theorem of Platis are the key ingredients in the construction of the
holomorphic extension of $\log\detp(\Delta)$.

\subsection{Holomorphic extension of $\log\detp(\Delta)$}
\label{S:holo ext}

The following is a key step in the proof of Theorem~\ref{main}.
\begin{proposition} \label{potential} Let $V$ and $W$ be domains in
  the complex space $\C^n$ diffeomorphic to the open unit ball.
  Consider $V \times W\subset\C^n\times\C^n$, with holomorphic
  coordinates $(z, w)$, and let $\p_z = dz^i \, \p_{z_i}$ and $\p_w =
  dw^j\, \p_{w_j}$. Suppose $\Omega$ is a holomorphic closed 2-form on
  $V\times W$ which is locally written as
  \begin{equation*}
    \Omega = \sum_{i,j} \Omega_{ij} dz^i \wedge dw^j .
  \end{equation*}
  Then there is a holomorphic function $q$ on $V \times W$ such that
  $\p_z \p_w q = \Omega$.
\end{proposition}
\begin{proof}
  Choose smooth polar coordinates on $V$ and $W$, and denote the
  centers of these coordinate systems by $z_0$ and $w_0$ respectively.
  Denote the radial line in polar coordinates from $z_0$ to the point
  $z\in V$ by $\horizontal(z)$; similarly, denote the radial line in
  polar coordinates from $w_0$ to the point $w\in W$ by
  $\vertical(w)$. More generally, if $c$ is a smooth chain in $V$, let
  $\horizontal(c)$ denote the cone on $c$ with vertex $z_0$, and
  similarly if $c$ is a smooth chain in $W$, let $\vertical(c)$ denote
  the cone on $c$ with vertex $w_0$.

  Define $q(z,w)$ by the formula
  \begin{equation*}
    q(z,w) = \int_{\horizontal(z) \times \vertical(w)} \Omega .
  \end{equation*}
  Since the chain $\horizontal(z) \times \vertical(w)$ varies smoothly
  as $(z,w)$ varies, the function $q(z,w)$ is smooth. Observe that $q$
  is unchanged by isotopies of the coordinate systems on $V$ and $W$
  which fix the centers $z_0$ and $w_0$, and that $q$ vanishes on
  $V\times\{w_0\}$ and on $\{z_0\}\times W$.

  If $c$ is a differentiable curve in $W$ parametrized by the interval
  $[0,t]$, we have by Stokes's theorem
  \begin{equation*}
    q(z,c(t)) - q(z,c(0)) = \int_{\horizontal(z)\times c} \Omega +
    \int_{\{z\}\times \vertical(c)} \Omega
    - \int_{\{z_0\} \times\vertical(c)} \Omega -
    \int_{\horizontal(z) \times \vertical(c)} d\Omega .
  \end{equation*}
  The second and third terms on the right-hand side vanish, since
  $\Omega$ vanishes when restricted to the $2$-simplex
  $\{z\}\times\vertical(c)$, and the last term vanishes since
  $d\Omega=0$. Taking the limit $t\to0$, we see that
  \begin{equation} \label{d_w} \iota(0,c'(0)) dq(z,c(0)) = -
    \int_{\horizontal(z)\times c(0)} \iota(0,c'(0))\Omega .
  \end{equation}
  Since $\Omega$ is holomorphic along $\{z\} \times W$, it follows
  that $q$ is holomorphic along $\{z\}\times W$ as well. A similar
  argument shows that $q$ is holomorphic along $V\times\{w\}$;
  combining these two calculations, we see that $q$ is holomorphic on
  $V\times W$.

  We now calculate $\p_w\p_zq$. By \eqref{d_w},
  \begin{equation*}
    \p_w q(z,w)  = -\sum_{i=1}^n dw^i  \int_{\horizontal(z)\times\{w\}}
    \iota(\p_{w^i}) \Omega .
  \end{equation*}
  If $c$ is a differentiable curve in $V$, parametrized by the
  interval $[0,t]$, we have by Stokes's theorem
  \begin{equation} \label{d_z}
    (\p_wq)(c(t),w) - (\p_wq)(c(0),w) = \sum_{i=1}^n dw^i \left( -
    \int_{c\times\{w\}} \iota(\p_{w^i}) \Omega
    + \int_{\horizontal(c)\times\{w\}} d\iota(\p_{w^i})\Omega \right) .
  \end{equation}
  The second term on the right-hand side vanishes. Indeed,
  \begin{align*}
    d \iota(\p_{w^i}) \Omega &= - \p_{z^k} \Omega_{ji} dz^k \wedge
    dz^j - \p_{w^k} \Omega_{ji} dw^k \wedge dz^j \\
    &= - \sum_{j<k} \bigl( \p_{z^k} \Omega_{ji} - \p_{z^j} \Omega_{ki}
    \bigr) \, dz^k \wedge dz^j - \p_{w^k} \Omega_{ji} dw^k \wedge dz^j \\
    &= - \p_{w^k} \Omega_{ji} dw^k \wedge dz^j .
  \end{align*}
  Restricting to $\horizontal(c)\times\{w\}$, this differential form
  vanishes.

  Taking $t\to0$ in \eqref{d_z}, we see that
  \begin{equation*}
    \iota(c'(0),0)d (\p_wq)(c(0),w) = - \sum_{i=1}^n dw^i \,
    \iota(c'(0),0) \iota(\p_{w^i}) \Omega(c(0),w) ,
  \end{equation*}
  or in other words, $\p_z\p_wq=\Omega$.
\end{proof}

From Proposition~\ref{form of Omega}, we know that the holomorphic
2-form $\Omega$ of Theorem~\ref{Platis} satisfies the hypotheses of
Theorem~\ref{potential}. Restricted to the diagonal
$\Teich(S)=\{w=\oz\}\subset\QF(S)$, the differential equation in
Theorem~\ref{potential} for the holomorphic function $q$ on $\QF(S)$
becomes
\begin{equation*}
  \p \dbar q = i \omegaWP.
\end{equation*}
Thus, the proof of Theorem~\ref{potential} gives a method of
constructing a K\"ahler potential for the K\"ahler form $i\omegaWP$ on
the Teichm\"uller space, using the extended form $\Omega$ to
quasifuchsian space.

\begin{example*} (See p.\ 214 in \cite{IT}) When $S$ has genus $1$,
  the Teichm\"uller space $\Teich(S)$ may be identified with the upper
  half plane $\H$, and
  \begin{equation*}
    \omegaWP = -i (z-\oz)^{-2} \, dz \wedge d\oz .
  \end{equation*}
  One easily finds the K\"ahler potential $q(z)=\log(z-\oz)$. The
  method used in the proof of Theorem~\ref{potential}, applied to the
  2-form $\Omega = (z-w)^{-2} \, dz \wedge dw$, yields the holomorphic
  function
  \begin{equation*}
    q(z,w) = \log(z-w) - \log(z_0-w) - \log(z-w_0) + \log(z_0-w_0)
  \end{equation*}
  on the quasifuchsian space $\H\times\overline{\H}$.
\end{example*}

Using the holomorphic function $q$ on $\QF(S)$, we now construct the
holomorphic extension of $\log\detp(\Delta)$. The holomorphic function
\begin{equation*}
  \tilde{q}(z,w) = \frac{1}{2} \bigl( q(z,w) + \oq(\ow,\oz) \bigr)
\end{equation*}
on $\QF(S)$ restricts to a real function $\tilde{q}$ on the diagonal
such that
\begin{equation*}
  \p \dbar \tilde{q} = i \omegaWP .
\end{equation*}
\begin{theorem}
  \label{holo det in quasifuchsian}
  There exists a unique holomorphic extension of $\log\detp(\Delta)$
  to the quasifuchsian space $\QF(S)$. In coordinates $(z,w)$ on
  $\QF(S)\cong\Teich(S)\times\overline{\Teich(S)}$, this extension has
  the form
  \begin{equation*}
    \log\detp(\Delta)(z,w) =  -\frac{1}{6\pi} \tilde{q}(z,w) +
    \log\det( (\tau(z)-\otau(\ow))/2i ) + f(z) + \of(\ow) .
  \end{equation*}
\end{theorem}
\begin{proof}
  By Theorem~\ref{curvature}, the one-form
  \begin{equation*}
    \alpha = \p \bigl( \log \detp(\Delta)  + \frac{1}{6\pi} \tilde{q} -
    \log\det(\Im\tau) \bigr)
  \end{equation*}
  is holomorphic. Since $\Teich(S)$ is simply connected, it follows
  that there is a differentiable function $f$ such that
  \begin{equation*}
    df = \alpha .
  \end{equation*}
  Since $\dbar f=\alpha^{0,1}=0$, $f$ is seen to be holomorphic. The
  theorem is now proved by analytically extending each of the
  functions $\det(\Im\tau)$, $\tilde{q}$, $f$ and $\of$ in the
  holomorphic factorization
  \begin{equation*}
    \log\detp(\Delta) = \log\det(\Im\tau) + C_g\tilde{q} + f + \of
  \end{equation*}
  on $\Teich(S)$ to $\QF(S)$. The holomorphic extension of $\tilde{q}$
  is evident, since it is by construction the restriction of the
  holomorphic function $\tilde{q}$ on $\QF(S)$. The function $f$ is
  extended to $f(z)$, the function $\of$ to $\of(\ow)$, and the
  function $\det(\Im\tau)$ to
   \begin{equation*}
     \log \det( ( \tau(z) - \otau(\ow) )/2i ) .
   \end{equation*}
   (Note that the matrix $\tau(z)-\otau(\ow)$ is everywhere invertible
   on $\QF(S)$.) The uniqueness of the holomorphic extension of
   $\log\detp(\Delta)$ follows from Lemma \ref{identity theorem}.
\end{proof}

It would not be hard, using this theorem, to give an explicit lower
bound for the radius of convergence of the real analytic function
$\log\detp(\Delta)$ on $\Teich(S)$.

\section{Holomorphic Extensions of Laplacians and Their Determinants}
\label{Ch:Delta}

In this section, we prove Theorem~\ref{T:main}.  In
Section~\ref{S:holo-ext-lap}, we construct the family
$\{\Delta_{\mu,\nu} \}$, and show that it satisfies properties (1) and
(2) in Theorem~\ref{T:main}. In Section~\ref{S:det-a}, we show the
property (3) of Theorem~\ref{T:main}. In Section~\ref{S:estimates}, we
provides several necessary estimates on quasiconformal mappings.
Using the results of Section~\ref{S:estimates}, we prove in
Sections~\ref{S:spectral cut} and \ref{S:eigenvalue bound} the results
which are used in Section~\ref{S:det-a}.  From now on, we
denote by $\p$ and $\dbar$ the Cauchy-Riemann operators $\frac{1}{2}
(\p_x-i\p_y)$ and $\frac{1}{2}(\p_x+i\p_y)$, respectively.

\subsection{The holomorphic extension $\Delta_{\mu,\nu}$ of the
  Laplacian} \label{S:holo-ext-lap}

In this subsection, we construct the family $\{\Delta_{\mu,\nu}\}$ of
elliptic second order differential operators of Theorem~\ref{T:main},
and demonstrate properties (1) and (2).

Unless otherwisely stated, we restrict our domain to $\H$, and denote
by $\mu$ and $\nu$ smooth Beltrami differentials on $\H$ (that is,
smooth complex valued functions on $\H$ satisfying $\|\mu\|_\infty$,
$\|\nu\|_\infty <1$). By $\hat{\mu}$ we denote a Beltrami differential
on the lower half plane $\overline{\H}$ defined by
$\hat{\mu}(z)=\omu(\oz)$.  Denote by $\op_\mu$ the operator $\op-\mu
\p$, and by $\p_{\omu}$ the operator $\p-\omu \op$.

The following definition is due to Ahlfors and Bers.

\begin{definition} \label{D:fm, fmn}
  Given a pair $(\mu,\nu)$ of Beltrami differentials on $\H$, denote
  by $\fmn : \C \to \C$ the unique continuous normalized solution
  (i.e.\ fixing $0$, $1$ and $\infty$) of the Beltrami equation on $\C$,
  \begin{equation*}
    \begin{cases}
      \op_\mu \fmn &= 0 , \qquad \Im z>0 , \\
      \op_{\hat{\nu}} \fmn &= 0 , \qquad \Im z<0 .
    \end{cases}
  \end{equation*}
  Let $f^\mu=f_{\mu,\mu}$.
\end{definition}

We have the following result of Ahlfors and Bers \cite{AB}.
\begin{lemma} \label{L:diffeo fmn}
  $\fmn$ is a homeomorphism of the Riemann sphere
$\Chat = \mathbb{C} \cup \{\infty \}$. In particular,
 it is an open
  embedding of $\H$ into $\C$, and $\p\fmn$ is nowhere zero on $\H$.
\end{lemma}

By complex conjugation of the Beltrami differential equation in
Definition~\ref{D:fm, fmn}, we see that
\begin{equation} \label{reflection of fmn}
  \overline{f_{\nu,\mu}(z)} = f_{\mu,\nu}(\oz) .
\end{equation}
In particular, $\ofm(z) = \fm(\oz)$ and thus $\fm$ maps $\H$ onto
$\H$.  In fact, $\fmn$ maps $\H$ onto $\H$ if and only if $\fm|_\R
\equiv f^\nu|_\R$.

In our construction of $\Delta_{\mu,\nu}$, we use the result of
Ahlfors and Bers that the normalized solutions of the Beltrami
equations depend analytically on the Beltrami differentials. The
following lemma summarizes what we need (see \cite{AB}).
\begin{lemma} \label{L:fmn holo, antiholo}
  For each $z \in \mathbb{H}$, $f_{\mu,\nu} (z)$,
  $\overline{f_{\nu,\mu}} (z)$, $\p
  f_{\mu,\nu} (z)$ and $\p\overline{f_{\nu,\mu}}(z)$,
  depend holomorphically
  on $\mu$ and anti-holomorphically on $\nu$.
\end{lemma}

Now, we start with the following key calculation in our construction
of $\Delta_{\mu,\nu}$. By Lemma~\ref{L:diffeo fmn} and the inequality
$|\mu|<1$, the function
\begin{equation} \label{alpha}
  \alpha = \frac{1}{(1-|\mu|^2)\pfm}
\end{equation}
is bounded on $\H$.
\begin{proposition} \label{41}
  \begin{equation*}
    (\fm)^*\p\op = |\alpha|^2 \left( - \mu \,\p^2 +
      (1+ |\mu|^2) \, \op\p - \omu\,\op^2 + (\op_\mu \log \alpha) \,
      \p + (\p_{\omu} \log \oalpha)  \, \op \right) .
  \end{equation*}
\end{proposition}

One easily sees that when $\mu=0$, the above formula for
$(\fm)^*\p\op$ reduces to $\p\op$.

In the proof of Proposition \ref{41}, we denote $f^\mu$ by $f$, and
$(f^\mu)^{-1}$ by $h$. By the chain rule applied to the equations
$h\circ f=z$ and $\oh\circ f=\oz$, and the Beltrami equation, we see
that
\begin{equation} \label{42}
  \begin{pmatrix}
    \p h \circ f & \p \oh \circ f \\ \op h \circ f & \op \oh \circ f
  \end{pmatrix} =
  \begin{pmatrix}
    \p f & \p \of \\ \op f & \op \of
  \end{pmatrix}^{-1}
  = \begin{pmatrix}
    \alpha  & - \omu\,\alpha  \\ - \mu\,\oalpha & \oalpha
  \end{pmatrix} .
\end{equation}

By the chain rule applied to the equation $\p h \circ f = \alpha$, we
see that
\begin{equation*}
  \begin{pmatrix}
    \p f & \p \of \\ \op f & \op \of
  \end{pmatrix}
  \begin{pmatrix}
    \p^2 h \circ f \\ \op\p h \circ f
  \end{pmatrix}
  =
  \begin{pmatrix}
    \p\alpha \\ \op\alpha
  \end{pmatrix}.
\end{equation*}
Applying \eqref{42}, we see that
\begin{equation} \label{43}
  \op\p h \circ f = |\alpha|^2 \, \op_\mu \log \alpha
\end{equation}

\begin{proof}[Proof of Proposition \ref{41}]
  If $u$ is a $C^{\infty}$ function on $\H$, then
  \begin{equation*}
    (\fm)^* \p \op u  = (\p \op ( u \circ h ) ) \circ f.
  \end{equation*}
  We have
  \begin{align*}
    \p\op(u\circ h) &= \p \left( (\p u\circ h) \op h + (\op u\circ h)
    \op \oh \right) \\
    &= (\p^2 u\circ h) \p h \op h + (\op\p u\circ h) \p \oh \op h \\
    &+ (\p\op u\circ h) \p h \op \oh + (\op^2 u\circ h) \p \oh \op \oh \\
    &+ (\p u\circ h) \p\op h + (\op u\circ h) \p\op\oh .
  \end{align*}
  Composing on the right with $f$, we see that
  \begin{align*}
    \left( \p\op(u\circ h) \right) \circ f &= \p^2 u (\p h\circ f)(\op
    h\circ f) + \op\p u (\p\oh\circ f) (\op h\circ f) \\
    &+ \p\op u (\p h\circ f) (\op\oh\circ f) + \op^2 u (\p\oh\circ f)
    (\op\oh\circ f) \\
    &+ \p u (\p\op h\circ f) + \op u (\p\op\oh\circ f) .
  \end{align*}
  Applying \eqref{42} and \eqref{43}, the proposition follows.
\end{proof}

We wish to find an extension of $(\fm)^*\p\op$ which is holomorphic in
$\mu$. Because the formula for $(\fm)^*\p\op$ contains quantities such
as $|\pfm|^2$ and $|\mu|^2$, simply replacing $\fm$ by $f_{\mu,\nu}$
does not give a holomorphic extension of $\p\op$. Nor do other simple
extensions, such as $(\fm)^*\p\,(f^\nu)^*\op$. On the other hand,
replacing $\fm$, $\ofm$ and $\omu$ by $\fmn$, $\ofnm$, and $\onu$,
respectively we obtain by Lemma~\ref{L:fmn holo, antiholo} an operator
which depends holomorphically on $\mu$ and anti-holomorphically on
$\nu$.

\begin{definition} \label{D:def Deltamn}
  Given a pair of Beltrami differentials $(\mu,\nu)$, let
  \begin{equation*}
  \alpha_{\mu,\nu} = \frac{1}{(1-\mu\onu) \pfmn}
  \end{equation*}
  Define a second order differential operator $\Delta_{\mu,\nu}$ on
  functions on $\H$ by the formula
  \begin{equation*}
    \Delta_{\mu,\nu} = (\fmn - \ofnm )^2 (\p\op)_{\mu,\nu}
  \end{equation*}
  where
  \begin{equation*}
    (\p\op)_{\mu,\nu} = \alpha_{\mu,\nu}\overline{\alpha_{\nu,\mu}}
    \bigl( - \mu\,\p^2 + (1 + \mu\onu) \, \op\p - \onu\,\op^2
    + (\op_\mu \log \alpha_{\mu,\nu}) \, \p + (\p_{\onu} \log
    \overline{\alpha_{\nu,\mu}}) \, \op \bigr) .
  \end{equation*}
\end{definition}

The principal symbol of $\Delta_{\mu,\nu}$ in complex coordinates
$(z,\zeta)$ on the cotangent bundle $T^*\H$, where $\sigma(\op) = i
\zeta$, equals
\begin{equation*} \label{principal symbol of Deltamn}
  \symbol(\Delta_{\mu,\nu})(\zeta) = - (\fmn - \ofnm)^2 \,
  \alpha_{\mu,\nu} \overline{\alpha_{\nu,\mu}} (\zeta - \mu \ozeta ) (
  \ozeta - \onu \zeta )
\end{equation*}
\begin{lemma}
  The differential operator $\Delta_{\mu,\nu}$ is elliptic for any
  pair of Beltrami differentials $(\mu,\nu)$.
\end{lemma}
\begin{proof}
  By \eqref{reflection of fmn}, we have
  \begin{equation*}
    f_{\mu,\nu} (z) - \overline{f_{\nu,\mu}(z)} = \fmn(z) - \fmn(\oz)
  \end{equation*}
  which is nowhere vanishing on $\mathbb{H}$,
since $\fmn$ is a homeomorphism of $\C$.
  The functions $\pfmn$ and $\opofnm$ are nowhere vanishing on $\H$ by
  Lemma~\ref{L:diffeo fmn}. We also have the bounds
  $\|\mu(z)\|_\infty , \|\nu(z)\|_\infty < 1$, and the lemma follows.
\end{proof}

The following theorem is immediate.
\begin{theorem} \label{T:holo. form of laplacian}
  The elliptic family $\Delta_{\mu,\nu}$ is holomorphic in $\mu$ and
  anti-holomorphic in $\nu$, and coincides with $(\fm)^*\Delta$ when
  $\mu = \nu $.
\end{theorem}

The following proposition shows that $\Delta_{\mu,\nu}$ is the unique
such family of operators.
\begin{proposition} \label{P:applying identity theorem}
  Let $A_{\mu,\nu }$ be a family of operators on $C^\infty(\H)$
  holomorphic in $\mu$ and anti-holomorphic in $\nu$. If $A_{\mu,\mu}
  = 0 $ for all $\mu$, then $A_{\mu,\nu} = 0 $ for all $\mu,\nu$.
\end{proposition}
We need an elementary lemma.
\begin{lemma} \label{L:identity theorem}
  Let $\phi (s, t) $ be a function of complex variables $s,t$ which is
  holomorphic in $s$ and anti-holomorphic in $t$.  Suppose $\phi(s,s)
  = 0$ for all $s$. Then $\phi(s,t) = 0$ for all $s,t$.
\end{lemma}

\begin{proof}[Proof of Proposition~\ref{P:applying identity theorem}]
  Fix $\psi\in C^\infty(\H)$, $z\in\H$, and Beltrami differentials
  $\mu,\nu$. Let $s,t$ be complex parameters. Then $\phi(s,t) =
  (A_{(1-s)\mu+s\nu,(1-t)\mu+t\nu}\psi)(z)$ is holomorphic in $s$ and
  anti-holomorphic in $t$, and $\phi(s,s) = 0$ for all $s$.  By
  Lemma~\ref{L:identity theorem}, $\phi(s,t) = 0 $ for all $s, t$.
  This shows the proposition.
\end{proof}

Now fix a Riemann surface $X_0$ and the corresponding Fuchsian group
$G$ of the covering map $\H \rightarrow X_0$.  We show that the
restriction of the family $\{\Delta_{\mu,\nu}\}$ to $G$-invariant
Beltrami differentials $\mu,\nu\in M^G$ on $\H$ induces a family of
elliptic differential operators on $X_0$.

\begin{lemma} \label{L:invar Delta mu nu}
  If $\mu\in M^G$ and $g\in G$, $g^*(f^\mu)^* \Delta = (f^\mu)^*\Delta$.
\end{lemma}
\begin{proof}
  By the invariance of the hyperbolic metric $m_0$ on $\H$ under
  conformal mappings, and by the invariance of $\mu$ under $G$, it is
  clear that the pull-back metric $(f^\mu)^* m_0$ is invariant under
  $G$. So the Laplacian $(f^\mu)^* \Delta$ associated to the pull-back
  metric $(f^\mu)^* m_0$ is also invariant under $G$.
\end{proof}

\begin{proposition} \label{P:invariance}
  For every $ g \in G$, and for every $\mu,\nu\in M^G$,
  $g^*\Delta_{\mu,\nu} = \Delta_{\mu,\nu}$.
\end{proposition}
\begin{proof}
  Fix $ g \in G$. The family of operators $g^*\Delta_{\mu,\nu} -
  \Delta_{\mu,\nu}$ is holomorphic in $\mu$ and anti-holomorphic in
  $\nu$, and by Lemma~\ref{L:invar Delta mu nu}, it vanishes for
  $\mu=\nu$. Therefore, by Proposition~\ref{P:applying identity
    theorem}, $g^*\Delta_{\mu,\nu} - \Delta_{\mu,\nu} =0$ for all
  $\mu,\nu\in M^G$.
\end{proof}

By Proposition~\ref{P:invariance} and the identification of $M^G$ with
$M(X_0)$, we have
\begin{theorem} \label{T:second order on X}
  There is a unique family $\{ \Delta _{\mu,\nu} \mid \mu,\nu\in M(X_0)
  \}$ of elliptic second order differential operators on $X_0$ which
  satisfies properties (1) and (2) of Theorem~\ref{T:main}.
\end{theorem}

\subsection{Determinant of $\Delta_{\mu , \nu}$} \label{S:det-a}
In this section, we consider the determinant of
$\Delta_{\mu,\nu}$ and establish
the property $(3)$ in Theorem~\ref{T:main}.
To define the determinant of $\Delta_{\mu, \nu}$,
we use the method of using complex powers of elliptic operators
developed by Seeley \cite{Se1}, \cite{Se2}, although we follow Shubin
\cite{Sh} more closely.  (See also \cite{KV}.)

For the Fuchsian group $G$ of $X_0$,  let $\S$ be the closure of a fixed
fundamental domain of $G$. Let $\T$ be the neighborhood of $\S$
consisting of the union of all translates of $\S$ by elements of $G$
whose intersection with $\S$ is nonempty.

\begin{definition}
  Given $0< k < 1$ and $E>0$, let
    \begin{equation*}
      M^G_{k,E} = \{ \mu \in M^G \mid \|\mu\|_\infty \leq k , \|
      \mu \|_{C^{2}(\T)} \leq E \}
    \end{equation*}
  where the $C^2$-norm is defined using the flat metric on $\H$.
\end{definition}

The following theorems will be proved in Sections~\ref{S:spectral cut},
\ref{S:eigenvalue bound}.

\begin{theorem}
  \label{T:spectral cut}
  Given $0<k<1$, $E>0$ and $0<\theta_0<\pi$, there is $\epsilon>0$
  such that if $\mu,\nu\in M^G_{k,E}$ and
  $\|\mu - \nu \|_{C^1 (\T)} < \epsilon$, then
  \begin{equation*}
    | \arg (\symbol(\Delta_{\mu,\nu})) | < \theta_0 .
  \end{equation*}
\end{theorem}

\begin{theorem}
  \label{T:eigenvalue bound}
  There exists a constant $C> 0$ such that for every  $\mu$,
  $\nu\in M^G_{k,E}$ and  for any nonzero eigenvalue
  $\lambda$ of $\Delta_{\mu,\nu}$ on $X_0=\H/G$,
  \begin{equation*}
    |\lambda| \, \geq \,  C-O(\|\mu-\nu\|_{C^2(\T)}).
  \end{equation*}
\end{theorem}

Fix $0 < \theta_0 < \pi$. For the rest of section
denote $\Delta_{\mu, \nu}$  by $A$ and assume that $(\mu , \nu)$ belongs to
\begin{equation*}
  N_\epsilon =
  \{ (\mu , \nu) \mid \mu , \nu \in  M^G_{k,E} \
  and \ \| \mu -\nu\|_{C^2 (\T)} \leq \epsilon \}
\end{equation*}
where $\epsilon > 0$ will be determined in the following.

\subsubsection{Determinant of $\Delta_{\mu , \nu}$}
\label{SS:det A}
By Theorem~\ref{T:spectral cut}, we know
that for sufficiently small $\epsilon$ the principal symbol
$\symbol(A) (x, \zeta ) $ does
not take values in the closed conical sector
\begin{equation*}
  \Lambda =\{ \lambda : \theta_0 \leq \text{arg} \lambda \leq 2\pi -
  \theta_0 \}
\end{equation*}
in the spectral plane $\C$ for any $(x, \zeta ) \in T^* S \setminus
S$. This condition ensures that $\Spec(A)\cap\Lambda$ is finite, so
there is a closed sector $\Lambda_0 \subset \Lambda$ which has only
zero spectrum inside.

By Theorem~\ref{T:eigenvalue bound}, for sufficiently small $\epsilon
> 0 $, there is $\rho > 0$ such that
\begin{equation*}
  \Spec(A) \cap \{ z \, | \, |z| < \rho \} \subset \{0 \}.
\end{equation*}
Given $\exp(i\theta) \in \Lambda_0$, let $\Gamma_{(\theta)}$ be the
contour $\Gamma_{1,\theta}(\rho) \cup \Gamma_{0,\theta}(\rho) \cup
\Gamma_{2,\theta}(\rho)$, where
\begin{align*}
  \Gamma_{1,\theta}(\rho) &= \{ x \exp(i\theta) \mid x \ge \rho \} , \\
  \Gamma_{0,\theta}(\rho) &= \{ \rho \exp(i\phi) \mid \theta > \phi >
  \theta - \pi \} , \\
  \Gamma_{1,\theta}(\rho) &= \{ x \exp (i(\theta-\pi)) \mid \rho \le x
  \} .
\end{align*}

 Denote by
$R_\lambda$ the resolvent $(A-\lambda I)^{-1}$.
Then for $\Re s < 0$, define
  \begin{equation*}
    (A _s)_{(\theta)} = \frac{i}{2\pi} \int_{\Gamma_{(\theta)}}
    \lambda^{s} \, R_\lambda \, d\lambda .
  \end{equation*}
By the symbol calculus of \cite{Sh}, $A_s $ is trace class for $\Re s
< - 1$. In the following, we omit $\theta$ from the notation for
$(A_s)_{(\theta)}$ and $\Gamma_{(\theta)}$.

For $s\in\C$, define the modified complex power $A^{s,o}$ of $A$ by
\begin{equation*}
    A^{s,o} = A^k A_{s-k}
\end{equation*}
where $k$ is an integer chosen so that $\Re s-k<0$.
To see that this definition does not depend on the choice of $k$,
consider the operator
\begin{equation*} \label{P_0}
  P_0 = \frac{i}{2\pi} \int_{|\lambda|=\rho} R_\lambda \, d\lambda .
\end{equation*}
Observe that $P_0^2 = P_0$, $P_0 A_s=0$,
  and that $P_0$ commutes with
$A$, $A_s$ and $A^{s,o}$.
Then the well-definedness of $A^{s,o}$
follows  since
\begin{equation*}
  A^k A_{-k} = A_{-k} A^k = 1 - P_0 .
\end{equation*}

The modified complex power $A^{s,o}$ has group property:
\begin{equation*}
  A^{s,o} A^{w,o} = A^{s+w,o} .
\end{equation*}

Following the arguments in \cite{Sh} (pp.\ 94--106), we may show that
the kernel $A^{-s,o}(x,y)\,dy$ of $A^{-s,o}$ can be meromorphically
extended to all of $\C$, with simple poles contained in the set
\begin{equation*}
  \Bigl\{ \frac{2-j}{2} \Bigm| j\ge0 \Bigr\} \setminus \{ - j \mid j\ge0 \} .
\end{equation*}
It follows that the meromorphic function
\begin{equation*}
  \Tr(A^{-s,o}) = \int_M A^{-s,o}(x,x) \, dx
\end{equation*}
is regular at $s=0$.
\begin{definition} \label{D:determinant of A}
  $\detp(A) =
  \exp ( -\p_s |_{s=0} \Tr (A^{-s, o}))$
\end{definition}

As remarked by Kontsevich and Vishik in \cite{KV}, a change in the
choice of contour $\Gamma_{\theta}$ changes $\p_s |_{s=0} \Tr
A^{-s,o}$ by an element of $2 \pi i\Z$. After taking the exponential,
the determinant $\detp(A)$ is well-defined.

We summarize our discussion in the following theorem.
\begin{theorem} \label{T:def det}
  There exists $\epsilon>0$ such that $\detp(\Delta_{\mu,\nu})$ is
  defined on $N_\epsilon$.
\end{theorem}

\subsubsection{Holomorphy of $\detp(\Delta_{\mu,\nu})$}
\label{SS:variation-det-Delta-mu-nu}

Suppose $A$ belongs to a differentiable family of operators all of
which satisfy the above conditions for a fixed contour $\Gamma$.  Then
we have the following well-known variation formula for the
determinant, which can be proved by symbol calculus of the kernel of
complex powers as in \cite{Sh}.
\begin{equation} \label{diff. det}
  d \log \detp(A) = \p_s |_{s=0} \Tr
  ( s A^{-s-1,o} \, dA )
\end{equation}

In order to argue from \eqref{diff. det} that
$\detp(\Delta_{\mu,\nu})$ is holomorphic with respect to $\mu$ and
$\nu$, we must clarify one subtle point: the contour $\Gamma$ must be
chosen so that the spectrum of the operator $\Delta_{\mu,\nu}$
does not cross it as we perform the differentiation.

Fix $\mu_1 , \nu_1 \in M(X_0)$ and $\delta >0$.  For complex numbers
$|s|, |t| < \delta$, let
\begin{equation*}
  ( \mu_s , \nu_t ) = (\mu + s \mu_1 , \  \nu + t \nu_1 )
  \ \in N_\epsilon
\end{equation*}
and denote $\Delta_{\mu_s , \nu_t}$ by $A(s,t)$ and
$\Delta_{\mu_s , \nu_t }- \lambda$ for $\lambda \in \Lambda$
by $A_\lambda (s, t)$.

\begin{lemma}
  If $\delta$ is sufficiently small, there exists $R>0$ such that the
  resolvent $A_\lambda (s,t)^{-1}$ is bounded on
  \begin{equation*}
    \Lambda_R = \{ \lambda \in \Lambda \mid |\lambda | \ge R \} .
  \end{equation*}
\end{lemma}
\begin{proof}
  Consider a parametrix $B_\lambda (s ,t)$ of $A(s,t)$ and consider
  the equation
  \begin{equation*}
    B_\lambda (s, t) A_\lambda (s,t) = I + C_\lambda (s,t) ,
  \end{equation*}
  where $C_\lambda (s,t)$ is a smoothing operator such that
  \begin{equation*}
    ( 1 + |\lambda |) \| C_\lambda (s,t) \|
  \end{equation*}
  is bounded. (See \cite{Sh} pp.85--86.)  By continuity of the kernel
  of $C_\lambda (s,t)$ with respect to $s,t$, we see that $ \|
  C_\lambda (s,t) \| $ is uniformly bounded for $|s|, |t| < \delta $,
  when $\delta$ is sufficiently small, and from this the existence of
  $R$.
\end{proof}

The boundedness of the resolvent $A_\lambda (s,t)^{-1}$ is an open
condition; thus, if the operator $A(0,0)$ has no eigenvalues in the
bounded domain
\begin{equation*}
  \{ z \in \Lambda_0 \mid \rho < |z| < R \} ,
\end{equation*}
then $A(s,t)$ has no eigenvalues in this domain either, for
sufficiently small $\delta$. Recall that the only eigenvalue of
$A(s,t)$ inside the disk $\{ z \mid |z| < \rho \}$ is $0$, for
sufficiently small $\delta$.

In conclusion, for each $(\mu , \nu) \in N_\epsilon$ we can choose a
contour $\Gamma$ in such a way that the only eigenvalue of
$\Delta_{\mu_s,\nu_t}$ inside $\Gamma$ is zero, for any small
variation $(\mu_s , \nu_t )$ of $(\mu , \nu)$ in $N_\epsilon$.
Since the determinant is independent of the choice
of the contour, we have
\begin{theorem}
  \label{T:(3)}
  The function $\detp(\Delta_{\mu,\nu})$ is holomorphic in the region
  $N_\epsilon$, where $\epsilon$ is chosen as in Theorem~\ref{T:def
    det}.
\end{theorem}

The property (3) in Theorem~\ref{T:main} is a direct consequence of
this theorem.  Note that the flat Euclidean norm
$\|\cdot\|_{C^{2}(\T)}$ for $M^G$ and the hyperbolic norm
$\|\cdot\|_{C^2(X_0)}$ for $M(X_0)$ are equivalent since $\T$ is a
finite cover of compact $X_0$.

\section{Estimates for quasiconformal mappings} \label{S:estimates}

We start by reviewing some basic facts about quasiconformal mappings
due to Ahlfors and Bers \cite{AB}. Given $p>2$, let $C_p>1$ be the
constant associated to $p$ by Ahlfors and Bers (see p.\ 386,
\cite{AB}); note that
\begin{equation*}
  \lim_{p\searrow2} C_p = 1 .
\end{equation*}
Fix $0<k<1$, and choose $p>2$ such that $C_p<1/k$. We abbreviate
$L^p(\C)$ to $L^p$. Let $\mu$ and $\nu$ be complex valued functions in
$L^\infty(\C)$ with norm $\|\mu\|_\infty,\|\nu\|_\infty \le k$.

\begin{definition}
  \label{D:w^mu}
\cite{AB}
  The normalized solution $w^\mu: \C \to \C$ of the Beltrami equation
  $\op_\mu w^\mu = 0$ is the unique continuous
solution which fixes $0$, $1$, and
  $\infty$.
\end{definition}
It is known that the function $w^\mu$ is a homeomorphism of the
Riemann sphere $\Chat=\C \cup\{\infty\}$.
If $w^\nu = w^\rho \circ w^\mu $, then
\begin{equation*}
  \rho = \left( \frac{\nu-\mu}{1-\omu\nu}
    \frac{\p w^\mu}{\overline{\p w^\mu}} \right) \circ (w^\mu)^{-1} ,
\end{equation*}
 and $(w^\mu)^{-1} =
w^{\tilde{\mu}}$ where
\begin{equation}
  \label{tilde}
  \tilde{\mu} = \left( - \mu \frac{\p w^\mu}{\overline{\p w^\mu}} \right)
    \circ (w^\mu)^{-1} .
\end{equation}
Note that $\| \tilde{\mu} \|_\infty = \| \mu \|_\infty $.

Denote the spherical distance in the extended complex plane by
$[z_1,z_2]$. By Lemma 16 of \cite{AB}, there are positive constants
$\alpha(k)$ and $c(k)$ such that
\begin{equation*}
  [ w^\mu(z_1) , w^\mu(z_2) ] \le c(k) \, [z_1,z_2]^{\alpha(k)} .
\end{equation*}

Let $D_R=\{z\in\C\mid |z|\le R\}$ be the disk of radius $R$ in $\C$.
Since the spherical and Euclidean distances are equivalent in compact
domains, we see that if $z_1,z_2\in D_R$, then
\begin{equation} \label{AB34+36}
  |w^\mu(z_1) - w^\mu(z_2) | \le c(k,R) \, |z_1 - z_2 |^{\alpha(k)} .
\end{equation}
In particular, taking $z_0=0$, we see that
\begin{equation} \label{AB36}
  \| w^\mu \|_{L^\infty(B_R) } \le c(k,R) \, R^{\alpha(k)} .
\end{equation}

We also have the following lemma. (See p.\ 398 of \cite{AB}.)
\begin{lemma} \label{L:Ahlfors Bers lemma} If $\mu$ and $\nu$ are
  Beltrami differentials on $\Chat$ with $\|\mu\|_\infty ,
  \|\nu\|_\infty \le k$, then for all $z\in\Chat$,
  \begin{equation*}
    [w^\mu(z) , w^\nu (z) ] \le C(k) \, \| \mu - \nu \|_\infty.
  \end{equation*}
\end{lemma}

In particular,
\begin{equation} \label{Ahlfors Bers}
  \| w^\mu - w^\nu \|_{L^\infty(D_R)} \le C(k,R) \, \|\mu-\nu\|_\infty .
\end{equation}

We will need the following interior Schauder estimates for the
operators $\dbar_\mu$.
\begin{proposition}
  \label{P:w^mu elliptic bounds}
  Fix a bounded open domain $\Omega$ in $\C$, a relatively compact
  open subset $\Omega_1\ssubset\Omega$, a positive integer $n$, and
  real numbers $0<\delta<1$, $0<k<1$, and $E>0$. Let $\mu$ and $\nu$
  be Beltrami differentials on $\C$ satisfying $\|\mu\|_\infty$,
  $\|\nu\|_\infty \leq k$ and $\|\mu\|_{C^{n-1,\delta}(\Omega)}$,
  $\|\nu\|_{C^{n-1,\delta}(\Omega)}\leq E$. Then there is a positive
  constant $C$, depending only on the above data, such that
  $\|w^\mu\|_{C^{n,\delta}(\Omega_1)} \leq C$ and
  \begin{equation*}
    \| w^\mu - w^\nu \|_{C^{n,\delta}(\Omega_1)} \leq C \{
    \|\mu-\nu\|_{C^{n-1,\delta}(\Omega)} + \|\mu-\nu\|_\infty \} .
  \end{equation*}
\end{proposition}
\begin{proof}
  As long as $\|\mu\|_{L^\infty(\Omega)}$ is bounded by $k<1$, the
  operators $\dbar_\mu$ are uniformly elliptic on $\Omega$, and we
  have the uniform Schauder estimates
  \begin{equation*}
    \| w^\mu \|_{C^{n,\delta}(\Omega_1)} \leq
    C \|w^\mu\|_{C^0(\Omega)} ,
  \end{equation*}
  from which $\|w^\mu\|_{C^{n,\delta}(\Omega_1)} \leq C$ follows by
  \eqref{AB36}.

  Note that
  \begin{equation*}
    \op_\mu (w^\mu - w^\nu ) = (\mu - \nu )\p w^\nu
  \end{equation*}
  This implies the uniform Schauder estimates
  \begin{equation*}
    \| w^\mu - w^\nu \|_{C^{n,\delta}(\Omega_1)} \leq C \{ \|w^\mu -
    w^\nu \|_{C^0(\Omega)} + \| (\mu -\nu ) \p w^\nu
    \|_{C^{n-1,\delta}(\Omega)} \} ,
  \end{equation*}
  and applying \eqref{Ahlfors Bers}, the desired estimate on
  $\|w^\mu-w^\nu\|_{C^{n,\delta}(\Omega_1)}$ follows.
\end{proof}

The goal of the rest of  section is to verify the following theorem.
\begin{theorem}
  \label{T:lower bound}
  Let $\Omega_1 \ssubset \Omega \ssubset \C$. Suppose that
  $\|\mu\|_\infty \leq k$ and that $\|\p\mu\|_{L^p(\Omega)}<\infty$.
  Then the normalized solution $w^\mu$ of the Beltrami equation of
  $\mu$ satisfies
  \begin{equation*}
    \inf_{\Omega_1 } | \p w^\mu | \geq C  e^{-C\|\p \mu \|_{L^p (\Omega)}} .
  \end{equation*}
\end{theorem}

We will first consider the case where $\mu$ has compact support; we
imitate the proof of Lemma 7 in \cite{AB}. First, we recall some
results from \cite{AB} on the inhomogeneous Beltrami equation.
\begin{definition}
  \label{w^mu sigma}
  For $\sigma\in L^p$, let $w^{\mu,\sigma}: \C \rightarrow \C$ be the
  unique solution of the inhomogeneous Beltrami equation $\op_\mu w =
  \sigma$ such that $w(0)=0$ and $\p w\in L^p$.
\end{definition}

Two properties of $w^{\mu,\sigma}$ which we will need are
\begin{align}
  \label{AB14-0}
  \| \p w^{\mu,\sigma} \|_p &\le \frac{C_p\,\|\sigma\|_p}{1-kC_p}
  \intertext{and}
  \label{AB15-0}
  | w^{\mu,\sigma}(z_1) - w^{\mu,\sigma}(z_2) | &\le \frac{c_p
    \|\sigma\|_p}{1-kC_p} \, |z_1 - z_2 |^{1-2/p} .
\end{align}
(For the definition of the constant $c_p$, see p.\ 386 of \cite{AB}.)

\begin{lemma}
  \label{P:cpt supp bound}
  Suppose that $\|\mu\|_\infty \leq k$ and that $\p\mu\in L^p$. If
  $\mu$ has support in $D_R$, there is a constant $C$, depending only on
  $R$, such that
  \begin{equation*}
    \inf_{z \in \mathbb{C}}
    | \p w^\mu | \geq \frac{1}{1+k} e^{-C \|\p\mu\|_p} .
  \end{equation*}
\end{lemma}
\begin{proof}
  Let $\lambda = w^{\mu,\p\mu}$. By \eqref{AB14-0},
  \begin{align}
    \label{AB14}
    \|\op\lambda\|_p \le C \|\p\mu\|_p ,
  \end{align}
  while by \eqref{AB15-0},
  \begin{equation*}
    |\lambda(z_1 ) - \lambda(z_2 ) | \le C \, \| \p \mu \|_p |z_1 - z_2
    |^{1-2/p} .
  \end{equation*}
  Since $\lambda(0) = 0$,
  \begin{equation}
    \label{AB15-1}
    |\lambda(z)| \le C \, \|\p\mu\|_p |z|^{1-2/p} .
  \end{equation}
  In particular, when $|z|\le R+1$,
  \begin{equation*}
    \| \lambda(z) \| \le C (R+1)^{1-2/p} \|\p\mu\|_p .
  \end{equation*}

  If $R+1<|z|<r$, then since $\op\lambda(z)=0$ for $|z|>R$, Green's
  formula shows that
  \begin{equation*}
    \lambda(z) = \frac{1}{2\pi i} \int_{|\zeta|=r}
    \frac{\lambda(\zeta)}{\zeta-z} d\zeta
    + \frac{1}{2\pi i} \int\int_{D_R}
    \frac{\op\lambda(\zeta)}{\zeta-z} d\zeta d\ozeta .
  \end{equation*}
  Thus
  \begin{equation*}
    d\lambda(z) = \frac{1}{2\pi i}
    \int_{|\zeta|=r} \frac{\lambda(\zeta)}{(\zeta-z)^2} d\zeta
    + \frac{1}{2\pi i} \int\int_{D_R}
    \frac{\op\lambda(\zeta)}{(\zeta-z)^2} d\zeta d\ozeta .
  \end{equation*}
  By \eqref{AB15-1},
  \begin{equation*}
    \left| \frac{1}{2\pi i}
      \int_{|\zeta|=r} \frac{\lambda(\zeta)}{(\zeta-z)^2} d\zeta \right|
    \le \frac{1}{2\pi(r-|z|)^2} \int_{|\zeta|=r}
    |\lambda(\zeta)|\,|d\zeta| \le C \,\frac{r^{1-2/p}}{(r-|z|)^2}
    \|\p \mu\|_p.
  \end{equation*}
  By \eqref{AB14},
  \begin{align*}
    \left| \frac{1}{2\pi i} \int\int_{D_R}
      \frac{\op\lambda(\zeta)}{(\zeta-z)^2} d\zeta d\ozeta \right| &
    \le \frac{1}{2\pi(|z|-R)^2}
    \int\int_{D_R} |\op\lambda(\zeta)| \, |d\zeta d\ozeta | \\
    & \le C \frac{(\pi R^2)^{1-1/p}}{2\pi(|z|-R)^2}
    \|\p\mu\|_{L^p(\Omega)} .
  \end{align*}
  Taking $r\to\infty$, we see that
  \begin{align*}
    |d\lambda(z)| \le C \frac{\|\p\mu\|_p}{(|z|-R)^2} .
  \end{align*}
  Let $\zhat=z/|z|$. It follows that
  \begin{align*}
    |\lambda(z)| &\leq |\lambda((R+1)\zhat)| + \int_{R+1}^r
    |d\lambda(s\zhat)| \, ds \\
    &\leq C \|\p\mu\|_p \left( (R+1)^{1-1/p} + \int_1^\infty
      \frac{ds}{s^2} \right) .
  \end{align*}
  In summary, we see that
  \begin{equation}
    \label{lambda bound}
    \|\lambda\|_\infty \le C \|\p\mu\|_p .
  \end{equation}

  Let $\rho = e^\lambda$. Since $\op\rho=\p(\mu\rho)$, there exists a
  $C^1$ function $f$ such that $\p f= \rho$ and $\op f = \mu\rho$. As
  remarked in \cite{AB}, $f$ is a homeomorphism on the extended
  complex plane $\Chat$ and $f(\infty) = \infty$. Clearly, the
  normalized solution $w^\mu$ is
  \begin{equation*}
    w^\mu(z)  = \frac{f(z) - f(0)}{f(1) - f(0)} ,
  \end{equation*}
  hence
  \begin{equation*}
    |\p w^\mu |  \geq  \frac{e^{-|\lambda|}}{|f(1) - f(0)|} .
  \end{equation*}
  The numerator is bounded below by \eqref{lambda bound}, while the
  denominator is bounded above using the mean value theorem:
  \begin{align*}
    |f(1) - f(0)| & \leq \sup_D (|\p f| + |\op f| ) \\
    &\leq (1+k) \sup_D e^{|\lambda|} \leq \ (1+ k) e^{C\|\p\mu\|_p} .
    \qedhere
  \end{align*}
\end{proof}

\begin{proof}[Proof of Theorem~\ref{T:lower bound}]
  Choose an open set $\Omega'$ such that $\Omega_1 \ssubset \Omega'
  \ssubset \Omega$.  Let $\eta$ be a $C^\infty$ cut-off function which
  equals $1$ on $\Omega'$ and $0$ outside $\Omega$. Let
  $\psi=w^{\eta\mu} \circ (w^\mu)^{-1}$. Note that $\op\psi=0$ on
  $w^\mu[\Omega']$. Thus, on $\Omega'$,
  \begin{equation*}
    \p w^{\eta\mu } = (\p\psi\circ w^\mu) \p w^\mu .
  \end{equation*}
  It follows that
  \begin{equation}
    \label{w^mu frac}
    |\p w^\mu | = \frac{|\p w^{\eta\mu}|}{|\p\psi\circ w^\mu|} .
  \end{equation}
  We must bound this below on $\Omega_1$. The numerator is bounded
  below by Lemma \ref{P:cpt supp bound}.

  To get an upper bound for the denominator of \eqref{w^mu frac},
  note
  \begin{equation*}
    \sup_{\Omega_1} | \p\psi \circ w^\mu | = \sup_{w^\mu[\Omega_1]}
    |\p\psi | .
  \end{equation*}
  Let $r = \dist(w^\mu[\Omega_1], \C\setminus w^\mu[\Omega'] )$.
  Since $\psi$ is holomorphic on $w^\mu[\Omega']$, we see that for
  $z\in w^\mu[\Omega_1]$,
  \begin{equation*}
    \p\psi(z) = \frac{1}{2\pi i} \int_{|\zeta-z|=r}
    \frac{\psi(\zeta)}{(\zeta-z)^2} d\zeta
  \end{equation*}
  and therefore,
  \begin{equation*}
    \sup_{\Omega_1} | \p\psi\circ w^\mu | \leq r^{-1}
    \sup_{w^\mu[\Omega']} |\psi| .
  \end{equation*}
  But, by \eqref{AB36},
  \begin{equation*}
    \sup_{w^\mu[\Omega']} |\psi| = \sup_{w^\mu[\Omega']} |w^{\eta\mu}
    \circ (w^\mu)^{-1} | = \sup_{\Omega'} |w^{\eta\mu}| \le C .
  \end{equation*}
  It remains to bound $r$ below.

  Recall the definition \eqref{tilde} of the Beltrami differential
  $\tilde{\mu}$. By \eqref{AB34+36} and \eqref{AB36}, if
  $z_1\in w^\mu[\Omega_1]$ and $z_2\in\C\setminus w^\mu[\Omega']$,
  there is a constant $\alpha(k)>0$ such that
  \begin{equation*}
    | w^{\tilde{\mu}}(z_1) - w^{\tilde{\mu}}(z_2) | \le C
    |z_1-z_2|^{\alpha(k)} .
  \end{equation*}
  From this, we have
  \begin{equation*}
    \dist( \Omega_1 , \C\setminus\Omega' ) \le C \dist(
    w^\mu[\Omega_1] , \C \setminus w^\mu[\Omega'] )^\alpha .
  \end{equation*}
  So $r \geq C \dist(\Omega_1,\C\setminus\Omega)^{1/\alpha}$.
\end{proof}

\section{Proof of  Theorem \ref{T:spectral cut}} \label{S:spectral
  cut}

Recall that
\begin{equation*}
  \symbol(\Delta_{\mu,\nu}) = - (\fmn - \ofnm )^2 \,
  ((1-\mu\onu)^2 \pfmn \opofnm )^{-1} (\zeta - \mu \ozeta ) ( \ozeta -
  \onu \zeta ) .
\end{equation*}
By invariance of $\Delta_{\mu , \nu}$ under $G$ we only need to estimate
the argument of this symbol on $\S$, and for this  we will  use the
results of Section~\ref{S:estimates}.

\subsection{Angle estimates for $(\zeta-\mu\ozeta)(\ozeta-\onu\zeta)$}

Let $\varrho=\mu-\nu$. Then
\begin{align*}
  (\zeta - \mu\ozeta) (\ozeta-\onu\zeta) &= (1+\mu\onu) \zeta\ozeta -
  \mu\ozeta^2 - \onu\zeta^2 \\
  &= (1+|\nu|^2) |\zeta|^2 - (\nu\ozeta^2 + \onu \zeta^2) + \varrho \onu
  |\zeta|^2 - \varrho \ozeta^2 .
\end{align*}
But
\begin{equation*}
  (1 + |\nu|^2 ) |\zeta|^2 - ( \nu\ozeta^2 + \onu\zeta^2 ) \geq
  (1-|\nu|)^2 |\zeta|^2 \geq (1-k)^2 |\zeta|^2 ,
\end{equation*}
while on $\S$,
\begin{equation*}
  | \varrho\onu|\zeta|^2 - \varrho\ozeta^2 |
  \leq  2 \epsilon |\zeta|^2 .
\end{equation*}
Therefore, for sufficiently small $\epsilon$, we have on $\S$
\begin{equation*}
  | \arg ( (\zeta - \mu \ozeta) ( \ozeta - \onu \zeta ) ) | =
  O(\epsilon) .
\end{equation*}

\subsection{Angle estimates for $-(\fmn-\ofnm)^2$}

We estimate $\arg(-(\fmn-\ofnm)^2)$ by means of the decomposition
\begin{equation*}
  \fmn - \ofnm = (\fm - \ofm) + (\fmn - \fm) + (\fm- \ofnm) .
\end{equation*}

\begin{lemma}
  \label{L:min of Im f^mu}
  There is a constant $C=C(k,\S)>0$ such that if $\mu$ is a Beltrami
  differential such that $\|\mu\|_\infty\le k$, and $z\in\S$,
  \begin{equation*}
    \Im f^\mu(z) > C .
  \end{equation*}
\end{lemma}
\begin{proof}
  Let $F(K)$ be the family of $K$-quasiconformal mappings from the
  $\H$ to itself fixing $0$, $1$ and $\infty$. In particular,
  $f^\mu\in F(K)$, with
  \begin{equation*}
    K = \frac{1+k}{1-k} ,
  \end{equation*}
  By Theorem 2.1 of \cite{L}, $F(K)$ is normal on $\H$, that is, every
  sequence of elements of $F(K)$ contains a subsequence which is
  locally uniformly convergent in $\H$. Let $y$ be the infimum
  \begin{equation*}
    y = \inf_{f\in F(K),z\in\S} \Im f(z) .
  \end{equation*}
  Choose sequence $(f_n)\in F(K)$ and $(z_n)\in \S$ such that
  \begin{equation*}
    \lim_{n\to\infty} \Im f_n(z_n) \to y .
  \end{equation*}
  Since $\S$ is compact, by passing to a subsequence, we may assume
  that $(z_n)$ converges to a limit $z_\infty\in\S$. Since $F(K)$ is
  normal, there is a subsequence which is locally uniformly convergent
  in $\H$, with continuous limit $f_\infty$ such that $\Im
  f_\infty(z_\infty)=y$. By Theorem 2.2 of \cite{L}, $f_\infty$ is
  $K$-quasiconformal, hence injective. Thus, $y>0$, since
  $f_\infty(z_\infty)$ is in the interior of $D$.
\end{proof}

It follows that
\begin{equation*}
  \inf_\S |\fm - \ofm | \geq C .
\end{equation*}
By \eqref{Ahlfors Bers},
\begin{equation*}
  \| \fm - \ofnm \|_{L^\infty(\S)} + \| \fmn - \fm \|_{L^\infty(\S)}
  \leq C \| \varrho \|_\infty \leq C \epsilon .
\end{equation*}
Therefore, for sufficiently small $\epsilon$, we have on $\S$,
\begin{equation*}
  | \arg(-(\fmn-\ofnm)^2) | = O(\epsilon) .
\end{equation*}

\subsection{Angle estimates for $1-\mu\onu$}

We have
\begin{equation*}
  1-\mu\onu = 1 - |\nu|^2 - \varrho\onu .
\end{equation*}
Since $1-|\nu|^2\geq 1-k^2$ and $|\varrho\onu| \leq \epsilon k$, we see
that
\begin{equation*}
  | \arg(1-\mu\onu) | = O(\epsilon) .
\end{equation*}

\subsection{Angle estimates for $\pfmn\opfmn$} \label{SS:angle-est2}

To estimate the argument of
\begin{align*}
  \pfmn \opofnm & = ( \p \fm + \pfmn - \p \fm )
  ( \overline{\p \fm }  + \opofnm - \overline{\p \fm}) \\
  & = |\p \fm|^2 + \p \fm (\opofnm - \overline{\p \fm})
  + \overline{\p \fm }(\pfmn - \p \fm ) \\
  & + (\pfmn - \p \fm ) (\opofnm - \overline{\p \fm}) ,
\end{align*}
we need a lower bound for $|\p\fm|$ and upper bounds for $\pfm$,
$\pfmn-\p\fm$ and $\opofnm-\opfm$. Theorem~\ref{T:lower bound},
applied with $\Omega_1=\S$ and $\Omega=\T$, implies that
\begin{equation*}
  \inf_\S | \pfm | \geq C.
\end{equation*}
By Proposition~\ref{P:w^mu elliptic bounds}, we have the
estimates $\|\fm\|_{C^{1,\delta}(\S)} < C$,
\begin{equation*}
  \| \fmn - \fm \|_{C^{1,\delta }(\S)} < C\| \mu - \nu \|_\infty ,
\end{equation*}
and $\| f_{\nu,\mu} - \fm \|_{C^{1, \delta}(\S)} < C ( \| \mu - \nu
\|_{C^1(\T)} + \| \mu - \nu \|_\infty )$. Therefore, for sufficiently
small $\epsilon$, we have on $\S$,
\begin{equation*}
    | \arg(\pfmn\opofnm) | = O(\epsilon) .
\end{equation*}
Combining the above estimates,
we obtain Theorem~\ref{T:spectral cut}.

\section{Proof of Theorem \ref{T:eigenvalue bound}}
\label{S:eigenvalue bound}

Before we proceed for the proof let us fix some notations. Let $m_0$
be the K\"ahler form of the standard hyperbolic metric on $\H$, let
\begin{equation*}
  m = (\fm)^* m_0 =  -2i \, \frac{|\p \fm |^2 (1- |\mu |^2)}{(\fm -
    \ofm )^2} dz \wedge d\oz
\end{equation*}
be the K\"ahler form of the pull-back hyperbolic metric by $f^\mu$
induced on $X_0$, and let $\Delta_m$ be the corresponding Laplacian.
Let $\langle - , - \rangle$ be the inner product on $L^2(X_0,m)$, and
let $\| \cdot \|_2$ be the $L^2$-norm. With respect to the frame $\{
dz, d\oz \}$ of the cotangent bundle $T^*X_0 \otimes \C$, the Hodge
star operator $\star$ (with respect to $m$) acts on $1$-forms as
\begin{equation*}
  \star
  \begin{pmatrix} a \\ b \end{pmatrix}
  = \frac{i}{(1-|\mu|^2)^2}
  \begin{pmatrix}
    2 \omu & - (1+ |\mu |^2 ) \\
    |\mu|^2+1 & -2\mu
  \end{pmatrix}
 \begin{pmatrix} \bar{a} \\ \bar{b} \end{pmatrix} .
\end{equation*}
From this, it is easy to see that for $u \in C^\infty (X_0)$,
\begin{align}\label{du comparison}
  \| du \|_2 ^2 = \int_{X_0} du \wedge \star du & \geq \int_{\S}
  (1-|\mu|^2 ) (|\p u|^2 + |\dbar u|^2) \, \frac{i \, dz \wedge d\oz}{2} \\
  & \geq (1-k) \| d u\|_{L^2 (X_0 , m_0 )}^2 . \notag
\end{align}
If $\nabla$ is the gradient operator of the metric $m$,
and $\nabla^*$ is its adjoint
with respect to $\langle - , - \rangle$, then $\Delta_m =
\nabla^*\nabla$; it follows that $\langle \Delta_mu , u \rangle = \|
\nabla u \|_2^2$.

\begin{lemma}
  \label{L:O estimates}
  Let $\epsilon = \| \mu - \nu \|_{C^{2}(\T)}$. Write
  $f=O_i(\epsilon^\ell)$ to denote that $f$ is a $C^i$ function (or
  tensor) such that $\| f \|_{C^i(\S)} \le C(k,E) \epsilon^\ell$. Then
  we have
  \begin{equation}
    \tag{1}
    \begin{aligned}
      \fmn &= O_2(1) , & (\p\fmn)^{-1} &= O_0(1) , & \fmn - \fm &=
      O_2(\epsilon) & \fmn - f^\nu &= O_2(\epsilon) , \\
      \alpha_{\mu,\nu} &= O_1(1) , & \dbar_\mu \log \alpha_{\mu,\nu}
      &= O_0(1) , & m / m_0 &= O_0 (1) ;
    \end{aligned}
  \end{equation}
  \begin{equation}
    \tag{2}
    \begin{aligned}
      \fmn - \ofnm & = \fm - \ofm + O_2(\epsilon) , &
      \alpha_{\mu , \nu} & = \alpha + O_1(\epsilon) \\
      \dbar_{\mu} \log \alpha_{\mu , \nu} & = \dbar_{\mu} \log \alpha
      + O_0(\epsilon) , & \p_{\onu} \log \overline{\alpha_{\nu , \mu}}
      & = \p_{\omu} \log \overline{\alpha} + O_0(\epsilon) ;
    \end{aligned}
  \end{equation}
  \begin{equation}
    \tag{3}
      \Delta_{\mu , \nu} = \Delta_m + O_1 (\epsilon) \dbar \p + O_1
      (\epsilon)\dbar^2 + O_0 (\epsilon) \p + O_0 (\epsilon) \dbar .
  \end{equation}
\end{lemma}
\begin{proof}
  (1) is by Proposition~\ref{P:w^mu elliptic bounds} and by
  Theorem~\ref{T:lower bound}, (2) is by straightforward calculation
  using (1), and (3) is by (1) and (2) and the definition of
  $\Delta_{\mu , \nu}$ (see Definition~\ref{D:def Deltamn}).
\end{proof}

From Lemma~\ref{L:O estimates} (3), we may write
\begin{equation*}
     \Delta_{\mu ,\nu} = \Delta_m + O_1(\epsilon) \nabla^2 +
   O_0(\epsilon) \nabla ,
\end{equation*}
where $O_i (\epsilon)$ is a tensor on $X_0$ whose $C^i$-norm is
bounded by $\epsilon$. Localization (by a partition of unity) and
integration by parts shows that
\begin{equation}\label{norm estimate Delta mu nu}
  \langle u , \Delta_{\mu ,\nu}u \rangle = \bigr( 1 + O(\epsilon)
  \bigl) \| \nabla u \|^2 + O(\epsilon) \| u\|_2^2 .
\end{equation}

Let $U$ be the space of constant functions on $X_0$, let $U^{\bot}$ be
its orthogonal complement in $L^2(X_0,m)$, and let
$\Delta_{\mu,\nu}^*$ be the adjoint of $\Delta_{\mu,\nu}$ with respect
to the metric $m$. If $f\in C^\infty(X_0)$, $\Delta_{\mu,\nu}^*f \in
U^\bot$.  Therefore, every eigenfunction $u$ of $\Delta_{\mu,\nu}^*$
with nonzero eigenvalue $\lambda$ belongs to $U^\bot$. If we let
\begin{equation*}
  v = u - \frac{\int_{X_0} u \, m_0}{\int_{X_0} m_0},
\end{equation*}
then clearly,
\begin{equation*}
  \frac{\| d v \|_2 ^2 }{\| v \|_2 ^2}
  \le \frac{ \| d u \|_2^2}{\|u\|_2^2 }
  = \frac{ \| \nabla u \|_2^2}{\|u\|_2^2 } .
\end{equation*}
By \eqref{du comparison},
\begin{align*}
  \| d v \|_{L^2 (X_0 , m_0 )}^2 \lesssim \| d v \|_2^2 .
\end{align*}
Since $m$ and $m_0$ are equivalent metrics, that is, $C^{-1}m_0 \le m
\le Cm_0$, we see that
\begin{align*}
  \| v \|_2 ^2 = \int_{X_0} |v|^2 m \sim \int_{X_0} |v|^2 \, m_0 .
\end{align*}
By the Poincar\'e inequality applied to $v$ for the metric $m_0$ on
$X_0$, we see that
\begin{align}
  \label{poincare type ineq}
  0 < C \le \frac{ \| \nabla u \|_2^2 }{\|u\|_2^2} ,
\end{align}
where the bound $C$ depends only on $k$ and $E$.

Since $u$ is an eigenfunction of $\Delta_{\mu,\nu}^*$ with nonzero
eigenvalue $\lambda$, we have by \eqref{norm estimate Delta mu nu},
\begin{align*}
  |\lambda| \, \| u \|_2 ^2 & = |\langle u , \Delta_{\mu,\nu}^* u
  \rangle | = | \langle \Delta_{\mu,\nu} u , u \rangle | \\
    & \geq \bigl( 1 - O(\epsilon) \bigr) \|\nabla u \|_2^2 -
  O(\epsilon ) \| u \|_2 ^2 .
\end{align*}
Therefore by the Poincar\'e inequality \eqref{poincare type ineq}, we
see that for sufficiently small $\epsilon$,
\begin{equation*}
  |\lambda | \geq C   - O(\epsilon)  .
\end{equation*}
This completes the proof of Theorem~\ref{T:eigenvalue bound}, since
the spectrum of $\Delta_{\mu,\nu}^*$ is the complex conjugate of the
spectrum of $\Delta_{\mu,\nu}$.

\bibliographystyle{plain}

\end{document}